\documentclass[final]{siamart220329}
\usepackage{amssymb}
\usepackage{amsmath}
\usepackage{amsfonts}
\usepackage{graphicx}
\usepackage{epstopdf}
\ifpdf
  \DeclareGraphicsExtensions{.eps,.pdf,.png,.jpg}
\else
  \DeclareGraphicsExtensions{.eps}
\fi
\usepackage{caption}

\usepackage[english]{babel}
\usepackage{hyphenat}

\newtheorem{assumption}[theorem]{Assumption}

\newtheorem{remark}[theorem]{Remark}
\numberwithin{equation}{section}

\newcommand{\rf}[1]{{\color{black}  #1}} 

\input{tcilatex}
\begin{document}

\title{Sampling from mixture distributions based on regime-switching
diffusions}

\author{M.V. Tretyakov\thanks{School of Mathematical Sciences, University of Nottingham,
Nottingham, UK (\email{Michael.Tretyakov@nottingham.ac.uk})
\funding{The author was partially supported by
Engineering and Physical Sciences Research Council (EPSRC) grant no. EP/X022617/1.}}}

\maketitle

\begin{abstract}
It is proposed to use stochastic differential equations with state-dependent
switching rates (SDEwS) for sampling from finite mixture distributions. An
Euler scheme with constant time step for SDEwS is considered. It is shown
that the scheme converges with order one in weak sense and also in the
ergodic limit. Numerical experiments illustrate the use of SDEwS for
sampling from mixture distributions and confirm the theoretical results.
\end{abstract}



\section{Introduction}

The main theme of this paper is to answer the question how to sample from a
finite mixture distribution using stochastic differential equations (SDEs).
Mixture distributions are probability distributions which can be written as
a superposition of component distributions. They are widely used in
statistical modelling and Bayesian inference arising, for instance, in
biology (e.g., population modelling, modelling DNA content, etc.),
engineering (e.g., modelling failures or defects in complex manufacturing
systems); image segmentation, clustering, financial application and many
others (see e.g. \cite{mix1,mix2,mix3}). Efficient sampling algorithms from
mixture distributions is of substantial practical importance while it is a
challenging problem.

For the last 30 years a lot of attention has been paid how to use SDEs for
efficient sampling from Gibbs distributions with densities proportional to $%
\exp (-U(x))$ with $x\in \mathbb{R}^{d}$ or $x\in D\subset \mathbb{R}^{d}$
or on Riemannian manifolds (see \cite{Tal90,RT96,MilTre07,MST10,MT,manifold}
and references therein). Despite applicable importance of mixture
distributions, an SDEs'-based approach has not been considered for sampling
from mixture distributions.

A natural SDEs' candidate for this sampling task is regime-switching
diffusions, in which the SDEs solution interacts with conditionally Markov
chain \cite{GN72,GN78b,Maobook,Yinbook}. As we show via a Fokker-Planck
equation argument (see Section~\ref{sec:prel}), to sample from a mixture
distribution, one needs to consider SDEs with state-dependent switching
rates (SDEwS). At the same time, most of the stochastic numerics literature
on SDEwS has dealt with constant switching rates only, except the works \cite%
{SDEwMS2,SDEwMS7,SDEwMS8,Yinbook}. As far as we are aware, in previous works
proofs of weak-sense convergence and convergence to ergodic limits for
numerical approximations for SDEwS were conducted predominately via first
proving mean-square (strong) convergence which then implies weak convergence
-- the route that has been typically avoided for usual SDEs, where instead
the link between SDEs and second-order linear parabolic and elliptic PDEs
was successfully used (see e.g. \cite%
{GN78a,GN85,Tal86,Tal90,MilTre07,MST10,MT}). In particular, the PDE route
allows to prove an optimal weak-sense order of convergence -- this approach
is exploited in this paper.

Although our main interest is in sampling from a finite mixture
distribution, we also consider finite-time weak-sense convergence for
approximations of SDEwS with state-dependent switching rates. To this end,
we recall that SDEwS find applications in mathematical modelling of
biological and ecological processes, in financial mathematics, economics,
wireless communications, physics, chemistry and engineering (see such
examples e.g. in \cite{Maobook,Yinbook,Lot24} and references therein), for
which weak-sense approximation is of practical interest. The proof of
finite-time weak convergence exploits the backward Kolmogorov equation
allowing to achieve an optimal weak-sense order of convergence.

In this paper we use an explicit Euler scheme with constant time stepping
for approximating SDEwS as a primer to illustrate the main concepts. It is
not difficult to consider other schemes and the proof technique (due to its
universality) for weak convergence presented here can be applied to them. We
consider the case of SDEwS with globally Lipschitz coefficients for
simplicity of the exposition. SDEwS with nonglobally Lipschitz coefficients
can be treated either by adapting the concept of rejecting exploding
trajectories from SDEs \cite{GNT04,MilTre07,MT} to SDEwS or by exploiting
implicit schemes (see e.g. \cite{SDEwMS4}) or tamed/balanced schemes as in
the case of usual SDEs (see e.g. \cite{Handy13,HJ15,MT}), which could be a
topic for a future work.

Structure of the rest of the paper is as follows. In Section~\ref{sec:prel}
we introduce SDEwS and the corresponding backward and forward Kolmogorov
equations. Then we derive SDEwS suitable for sampling from a given mixture
distribution. In Section~\ref{sec:meth} we write a numerical method and
state its properties: weak order of convergence is equal to one and in the
ergodic case convergence to ergodic limits in time step is also of order
one. In Section~\ref{sec:exp} we present results of numerical experiments.
Proofs of the convergence results are given in Section~\ref{sec:proof}.

\section{SDEs with switching\label{sec:prel}}

Let $\mathcal{M}=\{1,\ldots ,m_{0}\}$ and introduce the matrix
\begin{equation}
Q(x)=\left[
\begin{array}{cccc}
-q_{1}(x) & q_{12}(x) & \cdots & q_{1m_{0}}(x) \\
q_{21}(x) & -q_{2}(x) & \cdots & q_{2m_{0}}(x) \\
\vdots & \vdots & \cdot & \vdots \\
q_{m_{0}1}(x) & q_{m_{0}2}(x) & \cdots & -q_{m_{0}}(x)%
\end{array}%
\right] ,\ x\in \mathbb{R}^{d},\   \label{eq:Q}
\end{equation}%
where $q_{i}(x)$ and $q_{ij}(x)$ are continuous non-negative functions
bounded in $\mathbb{R}^{d}$ and the following relations hold:
\begin{equation}
\sum_{j\neq i}q_{ij}(x)=q_{i}(x),\ \ i\in \mathcal{M}.  \label{eq:qd}
\end{equation}

Consider the SDEwS%
\begin{equation}
dX=a(t,X(t),\mu (t))dt+\sigma (t,X(t),\mu (t))dw,\ \ X(t_{0})=x,\ \mu
(t_{0})=m,  \label{eq:sde}
\end{equation}%
where $X$ and $a$ are $d$-dimensional column-vectors, $\sigma $ is a $%
d\times d$-matrix, $w(t)=(w_{1}(t),\allowbreak \ldots ,w_{d}(t))^{\top }$ is a $d$%
-dimensional standard Wiener process, and $\mu (t)$ is a one-dimensional
continuous-time conditionally Markov process with right-continuous sample
paths taking values in $\mathcal{M}$ and
\begin{eqnarray}
P(\mu (t+\delta ) &=&i|\mu (t)=j,X(t)=x)=q_{ij}(x)\delta +o(\delta ),\ i\neq
j,  \label{eq:sde2} \\
P(\mu (t+\delta ) &=&i|\mu (t)=i,X(t)=x)=1-q_{i}(x)\delta +o(\delta )  \notag
\end{eqnarray}%
provided $\delta \downarrow 0.$ The Wiener process $w(t)$ and the
exponential clock governing transitions of $\mu (t)$ are independent.

We note that, using Skorokhod's representation \rf{(see e.g. \cite[p. 104]{Skoroh} and \cite[p. 29]{Yinbook} and \cite{SDEwMS8})}, the
process $\mu (t)$ can be expressed in terms of the Poisson random measure in
the following way. First, introduce the sequence of consecutive,
left-closed, right-open intervals $\Gamma _{ij}(x)$ for each $x\in \mathbb{R}%
^{d}:$ $\Gamma _{ii}(x)=\varnothing $ and for $i\neq j$%
\begin{eqnarray*}
\Gamma _{ij}(x) &=&\varnothing \ \text{if }q_{ij}(x)=0, \\
\Gamma _{ij}(x) &=&\left[ \sum_{l=i}^{i-1}q_{l}(x)+\sum_{l=1,l\neq
i}^{j-1}q_{il}(x),\sum_{l=i}^{i-1}q_{l}(x)+\sum_{l=1,l\neq
i}^{j}q_{il}(x)\right) .
\end{eqnarray*}%
Each non empty $\Gamma _{ii}(x)$ has length $q_{ij}(x).$ Next, \rf{define} the
function $F$: $\mathbb{R}^{d}\times \mathcal{M\times }\mathbb{R\rightarrow Z}
$%
\begin{equation*}
F(x,i,z):=\sum_{j\in \mathcal{M}}(j-i)\mathbb{I}_{\Gamma _{ij}(x)}(z).
\end{equation*}%
Then the process $\mu (t)$ satisfies the SDE
\begin{equation}
d\mu (t)=\int_{[0,L]}F(X(t),\mu (t-),z)N(dt,dz),  \label{eq:Poisrm}
\end{equation}%
where $L=m_{0}(\rf{m_{0}-1})\ell $ with $\ell =\max_{i\in \mathcal{M}}\sup_{x\in
\mathbb{R}^{d}}$ $q_{i}(x)$ and
$N(dt,dz)$ is a Poisson random measure with
intensity $dt\times \upsilon (dz)$ and $\upsilon (dz)$ being the Lebesgue
measure on $[0,L].$ The Poisson random measure $N(\cdot ,\cdot )$ is
independent of the Wiener process $w(\cdot ).$

The solution of (\ref{eq:sde}), (\ref{eq:sde2}) is a Markov process $(X,\mu)$
with the infinitesimal generator acting on functions $g\in C^{2}(\mathbb{R%
}^{d}\times \mathcal{M})$ \rf{(see \cite{GN78b} and \cite[p. 28]{Yinbook})}:%
\begin{gather*}
\mathcal{L}g(x;m):=\frac{1}{2}\sum_{l,k=1}^{d}\sigma _{lk}(t,x,m)\sigma
_{lk}(t,x,m)\frac{\partial ^{2}}{\partial x^{l}\partial x^{k}}%
g(x;m)+\sum_{l=1}^{d}a^{l}(t,x,m)\frac{\partial }{\partial x^{l}}g(x;m) \\
-q_{m}(x)g(x;m)+\sum_{j\neq m}q_{mj}(x)g(x;j),\ \ m\in \mathcal{M}.
\end{gather*}%
The Cauchy problem for the system of linear parabolic PDEs (the backward
Kolmogorov equation)
\begin{eqnarray}
\frac{\partial }{\partial t}u(t,x;m)+\mathcal{L}u(t,x;m) &=&0,\ x\in \mathbb{%
R}^{d},\ t<T,  \label{eq:BKE} \\
u(T,x;m) &=&\varphi (x;m),\ x\in \mathbb{R}^{d},  \notag \\
m &\in &\mathcal{M},  \notag
\end{eqnarray}%
has the probabilistic representation \cite{GN78b} \rf{(see also \cite[p. 52]{Yinbook})}:%
\begin{equation}
u(t,x;m)=E\varphi (X_{t,x,m}(T);\mu _{t,x,m}(T)),  \label{eq:Kac}
\end{equation}%
where $(X_{t,x,m}(s),\mu _{t,x,m}(s)),$ $s\geq t,$ is a solution of (\ref%
{eq:sde}), (\ref{eq:sde2}) and $\varphi (x;m)$ are sufficiently smooth
functions with growth not faster than polynomial (we provide precise
assumptions for convergence theorems in the next section).

The adjoint operator $\mathcal{L}^{\ast }$ to $\mathcal{L}$ acting on
sufficiently smooth functions $\rho (x;m)$ is:
\begin{gather}
\mathcal{L}^{\ast }\rho (x;m):=\frac{1}{2}\sum_{l,k=1}^{d}\frac{\partial ^{2}%
}{\partial x^{l}\partial x^{k}}\left( \sigma _{lk}(x,m)\sigma _{lk}(x,m)\rho
(x;m)\right)   \label{eq:adj} \\
-\sum_{l=1}^{d}\frac{\partial }{\partial x^{l}}\left(
a^{l}(x,m)\rho (x;m)\right)
-q_{m}(x)\rho (x;m)+\sum_{j\neq m}q_{jm}(x)\rho (x;j),\ \ m\in \mathcal{M}.
\notag
\end{gather}%
Consequently, the Fokker-Planck equation (forward Kolmogorov equation) takes
the form
\begin{eqnarray}
\frac{\partial }{\partial t}\rho (t,x;m) &=&\mathcal{L}^{\ast }\rho
(t,x;m),\ x\in \mathbb{R}^{d},\ t>t_{0},  \label{eq:FP} \\
\rho (t_{0},x;m) &=&\delta (x-x_{0})\delta_{m,i},\ x\in \mathbb{R}^{d},  \,m \in \mathcal{M},  \notag
\end{eqnarray}%
with
\begin{equation*}
\sum_{m=1}^{m_{0}}\int_{\mathbb{R}^{d}}\rho (t,x;m)=1.
\end{equation*}%
Here, $\{\rho (t,x;j)\}_{j=1}^{m_{0}}$ gives the density for $%
(X_{t_{0},x_{0},i}(t);\mu _{t_{0},x_{0},i}(t)).$

\subsection{Ergodic case}

Let $U(x;m)\geq 0,$ $x\in \mathbb{R}^{d},$ $m\in \mathcal{M}$, be
sufficiently smooth functions. Consider a given finite mixture distribution
with sufficiently smooth density
\begin{equation}
\rho (x)=\frac{1}{\mathcal{Z}}\sum_{m=1}^{m_{0}}\rho (x;m),  \label{eq:mix}
\end{equation}%
where
\begin{equation}
\rho (x;m)=\alpha _{m}\exp (-U(x;m)),  \label{eq:Gib}
\end{equation}%
$\alpha _{m}>0$ are weights, and the normalisation constant $\mathcal{Z:}$%
\begin{equation*}
\mathcal{Z=}\rf{\sum_{m=1}^{m_{0}}\int_{\mathbb{R}^{d}}\rho (x;m)dx}.
\end{equation*}%
As highlighted in the Introduction, mixture distributions are commonly
employed in statistical modelling and Bayesian inference.

\begin{assumption}
\label{ass:erg}The elements $q_{ij}(x),$ $i\neq j,$ of the matrix $Q(x)$ are
positive H{\"{o}}lder continuous functions for all $x\in \mathbb{R}^{d}$ and
they are bounded in $\mathbb{R}^{d}.$ The functions $U(x;m),$ $m\in \mathcal{%
M},$ are twice differentiable and for all $x\in \mathbb{R}^{d}$ and $m\in
\mathcal{M}:$%
\begin{equation}
-\nabla U(x;m)\cdot x\leq -c_{1}|x|^{2}+c_{2},  \label{eq:condU}
\end{equation}%
where $c_{1}>0$ and $c_{2}\in \mathbb{R}.$
\end{assumption}

Introduce the SDEwS%
\begin{equation}
dX=-\frac{1}{2}\nabla U(X,\mu (t))dt+dw(t),\ \ X(0)=x,\ \mu (0)=m,
\label{eq:erg}
\end{equation}%
and
\begin{eqnarray}
P(\mu (t+\delta ) &=&i|\mu (t)=j,X(t)=x)=q_{ij}(x)\delta +o(\delta ),\ i\neq
j,  \label{eq:erg3} \\
P(\mu (t+\delta ) &=&i|\mu (t)=i,X(t)=x)=1-q_{i}(x)\delta +o(\delta )  \notag
\end{eqnarray}%
under $\delta \downarrow 0.$

We recall \rf{\cite[p. 117]{HAS}} that we say that a process $Z(t)=(X(t),\mu (t))$ is
ergodic if there exists a unique invariant measure $\pi $ of $Z$ and
independently of $x\in \mathbb{R}^{d}$ and $m\in \mathcal{M}$ there exists
the limit
\begin{equation}
\lim_{t\rightarrow \infty }E\varphi (X_{x,m}(t);\mu
_{x,m}(t))=\sum_{j=1}^{m_{0}}\int_{\mathbb{R}^{d}}\varphi (x;j)\,\pi
(dx,j):=\varphi ^{erg}  \label{PA31}
\end{equation}%
for any function $\varphi (x;m)$ with polynomial growth at infinity in $x$;
and we say that a process $Z(t)$ is exponentially ergodic if for any $x\in
\mathbb{R}^{d}$ and $m\in \mathcal{M}$ and any function $\varphi $ with a
polynomial growth in $x$
\begin{equation}
\left\vert E\varphi (X_{x,m}(t);\mu _{x,m}(t))-\varphi ^{erg}\right\vert
\leq C\left( 1+|x|^{\kappa }\right) e^{-\lambda t},\ \ t\geq 0,  \label{PA34}
\end{equation}%
where $C>0$ and $\lambda >0$ are some constants.

Under Assumption~\ref{ass:erg}, the SDEwS (\ref{eq:erg}), (\ref{eq:erg3}) is
exponentially ergodic with $X(t)$ having (in the ergodic limit) a smooth
density $\rho (x)=\sum_{m=1}^{m_{0}}\rho (x;m)$ \rf{(see e.g. \cite{Majda,CMH,XuZhu} and also \cite[Sections 4.2 and 4.5]{Yinbook})}, where $\rho
(x;m)$ satisfy the stationary Fokker-Planck equations
\begin{equation}
\mathcal{L}^{\ast }\rho (x;m)=0,\ \ m\in \mathcal{M}.  \label{eq:statFPE}
\end{equation}%
Here the adjoint operator $\mathcal{L}^{\ast }$ (cf. (\ref{eq:adj})) is%
\begin{equation*}
\mathcal{L}^{\ast }\rho (x;m)=\frac{1}{2}\Delta \rho (x;m)+\frac{1}{2}\nabla
\cdot \left( \rho (x;m)\nabla U(x;m)\right) -q_{m}(x)\rho (x;m)+\sum_{j\neq
m}q_{jm}(x)\rho (x;j).
\end{equation*}

Now let us consider how to use the SDEwS (\ref{eq:erg}), (\ref{eq:erg3}) for
sampling from the mixture distribution (\ref{eq:Gib}). It is not difficult
to check that the functions $\rho (x;m)$ from (\ref{eq:Gib}) satisfy (\ref%
{eq:statFPE}) if they solve the following system%
\begin{equation*}
-q_{m}(x)\rho (x;m)+\sum_{j\neq m}q_{jm}(x)\rho (x;j)=0,\ m\in \mathcal{M},
\end{equation*}%
or equivalently%
\begin{equation*}
\sum_{j\neq m}\left[ q_{jm}(x)\rho (x;j)-q_{mj}(x)\rho (x;m)\right] =0,\
m\in \mathcal{M},
\end{equation*}%
which is satisfied if for $j\neq m$%
\begin{equation}
\frac{q_{jm}(x)}{q_{mj}(x)}=\frac{\rho (x;m)}{\rho (x;j)}.
\label{eq:q_choice}
\end{equation}%
Thus, we arrived at the proposition.

\begin{proposition}
\label{prop:rates}Let Assumption~\ref{ass:erg} hold for the functions $%
U(x,m) $ and for the matrix $Q(x)$ satisfying (\ref{eq:q_choice}).\ Then the
solution \rf{$(X_{0,x,m}(t);\mu _{0,x,m}(t))$} of the SDEwS (\ref%
{eq:erg}), (\ref{eq:erg3}) with $Q(x)$ satisfying (\ref{eq:q_choice}) has
the unique invariant measure $\rho (x)$ from (\ref{eq:mix}), (\ref{eq:Gib}).
\end{proposition}

We remark that a choice of functions $q_{jm}(x)$ satisfying (\ref%
{eq:q_choice}) is not unique. For instance, we can choose
\begin{equation}
q_{jm}(x)=\rho (x;m)  \label{eq:choice1}
\end{equation}%
or
\begin{equation}
q_{jm}(x)=\frac{\rho (x;m)}{\beta _{m}\beta _{j}},  \label{eq:choice2}
\end{equation}%
for some $\beta _{m}>0,$ etc. We emphasise that as in all Markov chain-type
samplers we do not need to know the normalisation constant $\mathcal{Z}$ of
the targeted distribution, we only need to know relative weights of the
mixture components. We note that under Assumption~\ref{ass:erg} both choices
(\ref{eq:choice1}) and (\ref{eq:choice2}) guarantee boundedness and
positivity of $q_{jm}(x).$ Functions $q_{jm}(x)$ are related to how fast \rf{$%
(X_{0,x,m}(t);\mu _{0,x,m}(t))$} converges to its ergodic
limit. Finding an optimal choice of $Q(x)$ satisfying (\ref{eq:q_choice}) is
an interesting theoretical and practical question which we leave for future
study.

Proposition~\ref{prop:rates} ascertain that the SDEwS (\ref{eq:erg}), (\ref%
{eq:erg3}) with $Q(x)$ satisfying (\ref{eq:q_choice}) can be utilised for
sampling from mixture distributions (\ref{eq:mix}), (\ref{eq:Gib}). In the
next section we construct a Markov chain approximating the SDEwS (\ref%
{eq:erg}), (\ref{eq:erg3}) and show that it can approximate associated
ergodic limits.

We note that in the case of the SDEwS (\ref{eq:erg}), (\ref{eq:erg3}), the
backward Kolmogorov equation (\ref{eq:BKE}) takes the form
\begin{eqnarray}
\frac{\partial }{\partial t}u(t,x;m)+\mathcal{L}u(t,x;m) &=&0,\ x\in \mathbb{%
R}^{d},\ t<T,  \label{eq:BKE2} \\
u(T,x;m) &=&\varphi (x;m),\ x\in \mathbb{R}^{d},  \notag \\
m &\in &\mathcal{M},  \notag
\end{eqnarray}%
where%
\begin{equation*}
\mathcal{L}g(x;m):=\frac{1}{2}\Delta g(x;m)-\frac{1}{2}\nabla U(x;m)\cdot
\nabla g(x;m)-q_{m}(x)g(x;m)+\sum_{j\neq m}q_{mj}(x)g(x;j).
\end{equation*}

\section{Numerical method, estimators and their properties\label{sec:meth}}

In Section~\ref{sec:mft} we consider an Euler method for the SDEwS (\ref%
{eq:sde}), (\ref{eq:sde2}) and state a theorem on its weak convergence with
order one. In Section~\ref{sec:merg} we apply the Euler method to the\ SDEwS
(\ref{eq:erg}), (\ref{eq:erg3}) and state error estimates for ensemble
averaging and time averaging estimators of ergodic limits corresponding to a
mixture distribution.

Introduce the equidistant partition of the time interval $[0,T]$ into $N$
parts with the step $h=T/N$: $t_{0}<t_{1}<\cdots <t_{N}=T$, $\
t_{k+1}-t_{k}=h.$

We recall and adapt to SDEwS the following definitions \rf{\cite[Chap. 2]{MT}}.

\begin{definition}
\label{Dat01} A function $f(x)$ belongs to the class $\mathbf{F}$ if there
are constants $K>0,$ $\kappa >0$ such that for all $x\in \mathbb{R}^{d}$ the
following inequality holds:
\begin{equation}
|f(x)|\leq K(1+|x|^{\kappa })\,.  \label{Da04}
\end{equation}%
If a function $f(x;s)$ depends not only on $x\in \mathbb{R}^{d}$ but also on
a parameter $s\in S,$ then $f(x;s)$ belongs to $\mathbf{F}$ $($with respect
to the variable $x)$ if an inequality of the type (\ref{Da04}) holds
uniformly in $s\in S.$
\end{definition}

\begin{definition}
\textit{If an approximation} $(\bar{X},\bar{\mu})$ of the solution $(X,\mu )$
to (\ref{eq:sde}), (\ref{eq:sde2}) \textit{is such that for }$p>0$
\begin{equation}
|Ef(\bar{X}(T);\bar{\mu}(T))-Ef(X(T);\mu (T))|\leq Ch^{p}  \label{A20}
\end{equation}%
\textit{for} $f\in \mathbf{F}$\textit{, then we say that the} weak order of
accuracy \textit{of the approximation} $(\bar{X},\bar{\mu})$ \textit{is} $p.$
The constant \rf{$C>0$} may depend on the coefficients of (\ref{eq:sde}), (\ref%
{eq:sde2}), on the function $f$ and on $T.$
\end{definition}

We introduce the following spaces: $C^{\frac{p+\epsilon }{2},p+\epsilon
}([0,T]\times \mathbb{R}^{d})$ (or $C^{p+\epsilon }(\mathbb{R}^{d})$) is a H{%
\"{o}}lder space containing functions $f(t,x)$ (or $f(x)$) whose partial
derivatives $\frac{\partial ^{i+|j|}f}{\partial t^{i}\partial
x^{j_{1}}\cdots \partial x^{j_{d}}}$ (or $\frac{\partial ^{|j|}f}{\partial
x^{j_{1}}\cdots \partial x^{j_{d}}}$) with $0\leq 2i+|j|<p+\epsilon $ (or $%
|j|<p+\epsilon $), $i\in \mathbb{N}\cup \{0\}$, $p\in \mathbb{N}\cup \{0\}$,
$0<\epsilon <1$, $j$ is a multi-index, are continuous in $[0,T]\times
\mathbb{R}^{d}$ (or $\mathbb{R}^{d}$) with finite H{\"{o}}lder norm $\mid
\cdot \mid ^{(p+\epsilon )}$ (see e.g. \rf{\cite[Chap. 1, Sec. 1]{LSU88}}). In what follows, for
brevity of the notation, we omit $\epsilon $ and write $C^{\frac{p}{2}%
,p}([0,T]\times \mathbb{R}^{d})$ (or $C^{p}(\mathbb{R}^{d})$) instead of $C^{%
\frac{p+\epsilon }{2},p+\epsilon }(\bar{Q})$ (or $C^{p+\epsilon }(\mathbb{R}%
^{d})$), which should not lead to any confusion.

\subsection{Finite-time case\label{sec:mft}}

We prove convergence of the numerical method considered in this section
under the following assumption on the coefficients of (\ref{eq:sde}), (\ref%
{eq:sde2}). At the same time, we note that the method itself can be used
when the assumption does not hold and studying this method under weaker
assumptions is a possible subject for a future work as discussed in the Introduction.

\begin{assumption}
\label{ass:ft}The coefficients $a(t,x,m)$ and $\sigma (t,x,m),$ $m\in
\mathcal{M},$ are $C^{1,2}\allowbreak([0,T]\allowbreak \times \mathbb{R}^{d})$ functions with
growth not faster than linear when $|x|\rightarrow \infty $; the elements $%
q_{ij}(x),$ $i\neq j,$ of the matrix $Q(x)$ are non-negative functions
bounded in $\mathbb{R}^{d}$ and belong to $C^{2}(\mathbb{R}^{d});$ and the
functions $\varphi (x;m),$ $m\in \mathcal{M},$ and their derivatives $\frac{%
\partial ^{|j|}\varphi }{\partial x^{j_{1}}\cdots \partial x^{j_{d}}},$ $%
|j|\leq 4,$ are continuous and belong to the class $\mathbf{F.}$
\end{assumption}

Under Assumption~\ref{ass:ft}, the solution $u(t,x;m)$ to the backward
Kolmogorov equation (\ref{eq:BKE}) satisfies the following estimate \rf{(see e.g.
\cite[Chap. 7, Sec. 10, Theorem 10.2]{LSU88} and \cite[Chap. 9, Sec. 6, Theorem 11]{Friedman})} for some $\kappa \geq 1:$
\begin{equation}
\sum\limits_{l=0}^{4}\sum\limits_{2i+|j|=l}\left\vert \frac{\partial ^{i+|j|}%
}{\partial t^{i} \partial x^{j_{1}} \cdots \partial x^{j_{d}}}%
u(t,x;m)\right\vert \leq K(1+|x|^{\kappa }),  \label{eq:FKest}
\end{equation}%
where $j$ is a multi-index and $K>0$ is independent of $x$ and $m$, it may
depend on $T.$

We also note that Assumption~\ref{ass:ft} is sufficient for existence of a
unique solution of (\ref{eq:sde}), (\ref{eq:sde2}) (see \rf{\cite[p. 30 and p. 34]{Yinbook} and} \cite{SDEwMS2}
for weaker conditions than Assumption~\ref{ass:ft} which are sufficient for existence of a
unique solution of (\ref{eq:sde}), (\ref{eq:sde2});  Assumption~\ref{ass:ft}
is used here for proving first-order weak convergence).

Consider the Euler scheme for the SDEwS (\ref{eq:sde}), (\ref{eq:sde2}):%
\begin{eqnarray}
X_{0} &=&x,\ \ \mu _{0}=m\in \mathcal{M},  \label{eq:Eul} \\
X_{k+1} &=&\rf{X_{k}}+ha(t_{k},X_{k},\mu _{k})+h^{1/2}\sigma (t_{k},X_{k},\mu
_{k})\xi _{k+1},  \notag \\
\mu _{k+1} &=&\left\{
\begin{array}{c}
j\neq \mu _{k}\ \ \text{with prob }q_{\mu _{k}j}(X_{k})h, \\
\mu _{k}\ \ \ \text{with prob }1-q_{\mu _{k}}(X_{k})h,%
\end{array}%
\right.  \notag
\end{eqnarray}%
where $\xi _{k+1}=(\xi _{k+1}^{1},\ldots ,\xi _{k+1}^{d})^{\top }$ and $\xi
_{k+1}^{i}$ are i.i.d. random variables with bounded moments and
\begin{equation}
E\xi _{k+1}^{i}=0,\ \ E\left[ \left( \xi _{k+1}^{i}\right) ^{2}\right] =1,\
\ E\left[ \left( \xi _{k+1}^{i}\right) ^{3}\right] =0.  \label{eq:rn}
\end{equation}%
The conditions (\ref{eq:rn}) are, e.g. satisfied for normal distribution
with zero mean and unit variance or for discrete random variables taking
values $\pm 1$ with probability $1/2.$ The random variables $\xi _{k+1}^{i}$
and the random variable responsible for switching $\mu _{k+1}$ are mutually
independent.

The following finite-time weak convergence theorem is proved in Section~\ref%
{sec:ftconv}\ for the scheme (\ref{eq:Eul}).

\begin{theorem}
\label{thm:ft} Let Assumption~\ref{ass:ft} hold. Then the scheme (\ref%
{eq:Eul}) for the SDEwS (\ref{eq:sde}), (\ref{eq:sde2}) converges with weak
order one, i.e.,
\begin{equation*}
\left\vert E\varphi (X_{N};\mu _{N})-E\varphi (X_{t,x,m}(T);\mu
_{t,x,m}(T))\right\vert \leq Ch(1+|x|^{\kappa }),
\end{equation*}

where $C>0$ and $\kappa \geq 1$ are independent of $h$.
\end{theorem}

\begin{remark}
Since the proof of Theorem~\ref{thm:ft} utilises closeness of the
semi-groups associated to the SDEwS (\ref{eq:sde}), (\ref{eq:sde2}) and the
Euler scheme (\ref{eq:Eul}), it is essentially universal \rf{like} all weak-sense
convergence proofs of this type (see e.g. \rf{\cite[Chap. 2]{MT}}) and can be extended to
other weak schemes for (\ref{eq:sde}), (\ref{eq:sde2}). E.g., it is of
interest to find a way of constructing second-order weak schemes for SDEs
with switching.
\end{remark}

\subsection{Ergodic case\label{sec:merg}}

Now we consider how to approximate ergodic limits associated with the SDEwS (%
\ref{eq:erg}), (\ref{eq:erg3}). To prove convergence of the considered
numerical method to ergodic limit with first order in $h$, we need a
stronger assumption than Assumption~\ref{ass:erg}.

\begin{assumption}
\label{ass:erg2} The functions $U(x,m)\in C^{\infty }(\mathbb{R}^{d}),$ $%
m\in \mathcal{M},$ and satisfy (\ref{eq:condU}) and $\nabla U(x,m)$ have
growth not faster than linear when $|x|\rightarrow \infty $; the elements $%
q_{ij}(x)$ of the matrix $Q(x)$ are positive functions bounded in $\mathbb{R}%
^{d}$ and belong to $C^{\infty }(\mathbb{R}^{d});$ and the functions $%
\varphi (x;m)\in C^{\infty }(\mathbb{R}^{d}),$ $m\in \mathcal{M},$ and they
and their derivatives belong to the class $\mathbf{F.}$
\end{assumption}

Similarly how it is proved in \cite{Tal90} for the case of the usual ergodic
SDEs, there are constants $K>0$, $\lambda >0$ and $\kappa \geq 1$
independent of $T,$ $x$ and $m$ such that the following bound holds for the
solution $u(t,x;m)$ to the backward Kolmogorov equation (\ref{eq:BKE2}):%
\begin{equation}
\sum\limits_{l=1}^{4}\sum\limits_{2i+|j|=l}\left\vert \frac{\partial ^{i+|j|}%
}{\partial t^{i}\partial x^{j_{1}} \cdots \partial x^{j_{d}}}%
u(t,x;m)\right\vert \leq K(1+|x|^{\kappa })e^{-\lambda (T-t)},
\label{eq:FKerg}
\end{equation}%
where $j$ is a multi-index. We note that under Assumption~\ref{ass:erg2} the
above estimate is true for any order of derivatives of the solution $%
u(t,x;m) $ \cite{Tal90}; we write (\ref{eq:FKerg}) for $2i+|j|\leq 4$
because it is sufficient for proofs.

In the case of the SDEwS (\ref{eq:erg}), (\ref{eq:erg3}), the Euler scheme (%
\ref{eq:Eul}) takes the form
\begin{eqnarray}
X_{0} &=&x,\ \ \mu _{0}=m\in \mathcal{M},  \label{eq:hmeth} \\
X_{k+1} &=&\rf{X_{k}}-\frac{h}{2}\nabla U(X_{k};\mu _{k})+h^{1/2}\xi _{k+1},
\notag \\
\mu _{k+1} &=&\left\{
\begin{array}{c}
j\neq \mu _{k}\ \ \text{with prob }q_{\mu _{k}j}(X_{k})h, \\
\mu _{k}\ \ \ \text{with prob }1-q_{\mu _{k}}(X_{k})h,%
\end{array}%
\right.  \notag
\end{eqnarray}%
where $\xi _{k+1}$ satisfy (\ref{eq:rn}).

Ergodic limits can be computed using ensemble average or time averaging \cite%
{Tal90,MilTre07,MT}. Ensemble averaging is based on proximity of $E\varphi
(X_{0,x,m}(T);\mu _{0,x,m}(T))$ to $\varphi ^{erg}$ for a sufficiently large
$T$ (see (\ref{PA34})). In practice the expectation $E\varphi (X_{N};\mu
_{N})$ approximating $E\varphi (X_{0,x,m}(T);\mu _{0,x,m}(T))$ is realised
using the Monte Carlo technique which can benefit from parallel computations
on GPUs and multiple CPUs. For faster exploration of the phase space and
hence for faster convergence to the ergodic limit, the initial points $x$
and $m$ of the process $(X_{0,x,m}(t),\mu _{0,x,m}(t))$ can be taken random
using a suitable distribution, e.g. $m$ can be distributed uniformly (see
further discussion in \cite{MilTre07} in the case of usual SDEs and here in
Section~\ref{sec:exp}).

The next theorem gives the error of approximating the ergodic limit using
ensemble-averaging, in other words, it gives bias of the ensemble-averaging
estimator $\hat{\varphi}^{erg}$ for \rf{the ergodic limit $\varphi ^{erg}$ defined in (\ref{PA31})}:
\begin{equation}
\hat{\varphi}^{erg}=\frac{1}{M}\sum_{k=1}^{M}\varphi (X_{N}^{(k)};\mu
_{N}^{(k)}),  \label{S1}
\end{equation}%
where $M$ is the number of independent realizations of $(X_{N},\mu _{N}).$
The error estimate is proved in Section~\ref{sec:ergproof}.

\begin{theorem}
\label{thm:erg1} Let Assumption~\ref{ass:erg2} hold. Then the scheme (\ref%
{eq:hmeth}) for the SDEwS (\ref{eq:erg}), (\ref{eq:erg3}) converges with
order one \rf{to the ergodic limit $\varphi ^{erg}$}, i.e.,
\begin{equation}\label{eq:thmerg1}
\left\vert E\varphi (X_{N};\mu _{N})-\varphi ^{erg}\right\vert \leq C\left[
h+(1+|x|^{\kappa })\exp (-\lambda T)\right] ,
\end{equation}
where $C>0$, \rf{$\kappa \ge 1$} and $\lambda >0$ are independent of $h$, \rf{$T$ and $x$; $C$ and $\kappa$ depend on the choice of $\varphi(x;m)$}.
\end{theorem}

\rf{We note that this theorem and its proof hold if instead of positivity of all elements of the matrix $Q(x)$
(see Assumption~\ref{ass:erg2}), we assume
that they are nonnegative and $Q(x)$ is irreducible for every $x \in \mathbb{R}^{d}$.
We keep the positivity assumption in the view of Proposition~\ref{prop:rates} and that the key objective of this paper is
to address the problem of sampling from mixture distributions.
}

Due to ergodicity, we have \rf{(see e.g. \cite[p. 110]{HAS})}
\begin{equation}
\lim_{t\rightarrow \infty }\frac{1}{t}\int\limits_{0}^{t}\varphi
(X_{0,x,m}(s);\mu _{0,x,m}(s))ds=\varphi ^{erg}\ \ a.s.  \label{PB51}
\end{equation}%
Then for a sufficiently large $\tilde{T}$:%
\begin{equation}
\varphi ^{erg}\approx \frac{1}{\tilde{T}}\int\limits_{0}^{\tilde{T}}\varphi
(X_{0,x,m}(s);\mu _{0,x,m}(s))ds\approx \check{\varphi}^{erg}:=\frac{1}{L}%
\sum_{l=1}^{L}\varphi (X_{l};\mu _{l}).  \label{PB52}
\end{equation}%
It follows from (\ref{PA34}) that the bias of the estimator $\frac{1}{\tilde{%
T}}\int\limits_{0}^{\tilde{T}}\varphi (X_{0,x,m}(s);\mu _{0,x,m}(s))ds$ is
of order $O(1/\tilde{T}):$
\begin{equation}
\left\vert E\left( \frac{1}{\tilde{T}}\int\limits_{0}^{\tilde{T}}\varphi
(X_{0,x,m}(s);\mu _{0,x,m}(s))ds\right) -\varphi ^{erg}\right\vert \leq
C\left( 1+|x|^{\kappa }\right) \frac{1}{\tilde{T}},  \label{eq:biasTime}
\end{equation}%
where $C>0$ is independent of $\tilde{T}$, $x$ and $m.$

\begin{theorem}
\label{thm:erg2} Let Assumption~\ref{ass:erg2} hold. Then the bias of the
estimator $\check{\varphi}^{erg}$ has the error estimate%
\begin{equation*}
\left\vert E\check{\varphi}^{erg}-\varphi ^{erg}\right\vert \leq C\left( h+%
\frac{1+|x|^{\kappa }}{\tilde{T}}\right) ,
\end{equation*}%
where $C>0$ \rf{ and $\kappa \ge 1$ are} independent of $h$, \rf{$x$}, and $\tilde{T}$; \rf{$C$ and $\kappa$ depend on the choice of $\varphi(x;m)$}.
\end{theorem}

Both estimators $\hat{\varphi}^{erg}$ and $\check{\varphi}^{erg}$ have the
three errors: numerical integration error controlled by the time step $h$
(and in general by a choice of numerical method too); the variance which is
controlled by the number of independent realisations $M$ in the case of the
ensemble-averaging estimator $\hat{\varphi}^{erg}$ and by the time $\tilde{T}%
=Lh$ in the case of $\check{\varphi}^{erg};$ and error due to proximity to
the stationary distribution which is controlled by the integration time $T=Nh
$ for $\hat{\varphi}^{erg}$ and by $\tilde{T}=Lh$ for $\check{\varphi}^{erg}.
$

\rf{
\begin{remark}
Let us briefly compare the two estimators.
First, it is not difficult to establish that  \begin{equation*}
\left\vert E\check{\varphi}^{erg}-E\hat \varphi ^{erg}\right\vert \leq C\left( h+%
(1+|x|^{\kappa })\left[\frac{1}{\tilde{T}}+\exp(-\lambda T)\right]\right) ,
\end{equation*}%
where $C>0$ is independent of $h$, \rf{$x$}, $T$, and $\tilde{T}$.
Next, the computational cost of $\hat \varphi ^{erg}$ is proportional to $NM$ and of $\check \varphi ^{erg}$ $\propto$ $L$.
Since the variance (which is just the Monte Carlo error) of $\hat \varphi ^{erg}$ is ${\cal O}(1/M)$,
 the mean-square error of
$\hat \varphi ^{erg}$ is (see, e.g. \cite[Sec. 2.5.2]{MT})
$\sqrt{E [ \hat \varphi ^{erg}-\varphi ^{erg}] ^{2}}={\cal O}(1/\sqrt{M}+h+\exp(-\lambda Nh))$, where we used Theorem~\ref{thm:erg1}.
Since the variance of $\check \varphi ^{erg}$ is ${\cal O}(1/Lh)$ (see \cite{MST10} and Remark~\ref{rem:Poisson}),
the mean-square error of $\check \varphi ^{erg}$ is ${\cal O}(1/\sqrt{Lh}+h)$, where we used Theorem~\ref{thm:erg2}.
Consequently, if we fix a tolerance level $\varepsilon > 0$ for the mean-square errors, then
the cost of reaching this level is ${\cal O}(\varepsilon ^{-3} \log(1/\varepsilon))$ for $\hat \varphi ^{erg}$
and ${\cal O}(\varepsilon ^{-3})$
for $\check \varphi ^{erg}$
if we optimally choose $h \propto \varepsilon$, $M \propto \varepsilon ^{-2}$, $N \propto \varepsilon ^{-1} \log(1/\varepsilon)$ and $L \propto \varepsilon ^{-3}$.
Traditionally, in statistics and molecular dynamics, time-averaging estimators have been used.
However, with growing availability of servers with multiple CPUs and powerful GPU cards, the ensemble-averaging estimators
are becoming more attractive because they perfectly suit parallel computing.
We remark that the discussed complexity is asymptotic and constants in the error estimates
can play an important role for computational efficiency of the estimators.
From this angle and as discussed before, we note that the initial condition $(X_0, \mu _0)$ can be taken random to speed up convergence
of $\hat \varphi ^{erg}$ to $\varphi ^{erg}$, i.e. reduce the constant in front of $\exp(-\lambda T)$ in the error estimate (see e.g. such a discussion in \cite{MilTre07}).
\end{remark}
}

\section{Numerical experiments\label{sec:exp}}

We conduct experiments with the ensemble-aver\hyp{}aging estimator $\hat{\varphi}%
^{erg}$ which allows to fully benefit from parallel computing. The Monte
Carlo error presented for the ensemble-averaging estimator $\hat{\varphi}%
^{erg}$ in this section is given as usual by $2\sqrt{\bar{D}_{M}/M}$ with
\begin{equation}
\bar{D}_{M}=\frac{1}{M}\sum_{m=1}^{M}\left[ \varphi (X_{N}^{(k)};\mu
_{N}^{(k)})\right] ^{2}-\left[ \hat{\varphi}^{erg}\right] ^{2}.
\label{eq:MC}
\end{equation}

We test the first-order of weak convergence of the scheme (\ref{eq:hmeth})
on sampling from a univariate Gaussian mixture with two components, a
mixture of Gaussian and distribution with nonconvex potential, and a
multivariate Gaussian mixture. In all the experiments, we use the choice (%
\ref{eq:choice1}) and we uniformly distribute the initial state $\mu (0)$
and take $X(0)$ equal to the mean of the $\mu (0)$th component in the mixture. \medskip

\noindent \textbf{Example 1.} Consider a mixture of two univariate Gaussian
distributions:
\begin{equation}
\rho (x)=\frac{1}{\mathcal{Z}}\left[ \alpha _{1}\exp \left( -\frac{(x-%
\mathrm{m}_{1})^{2}}{2\sigma _{1}^{2}}\right) +\alpha _{2}\exp \left( -\frac{%
(x-\mathrm{m}_{2})^{2}}{2\sigma _{2}^{2}}\right) \right] ,\ \ x\in \mathbb{R}%
,  \label{eq:uGm2}
\end{equation}%
where $\mathcal{Z}$ is the normalisation constant. The parameters used are
given in Table~\ref{Tparam1} and for test purposes we take $\varphi
(x;m)=x^{2}$ for $m=1,2.$ The exact value $\varphi ^{erg}=4.87500\ (5\ $d.p.$%
).$

\begin{table}[h] \centering%
\caption{Example 1. The parameters of the mixture (\ref{eq:uGm2}).
\label{Tparam1}}\smallskip
\begin{tabular}{l|ll}
\hline
$m$ & $1$ & $2$ \\ \hline
$\alpha _{m}$ & $0.5$ & $0.4$ \\ \hline
$\mathrm{m}_{m}$ & $0$ & $3$ \\ \hline
$\sigma _{m}$ & $2$ & $0.5$ \\ \hline
\end{tabular}%
\end{table}%

\begin{figure}[tbp]
\begin{center}
\par
\includegraphics[width=8cm]{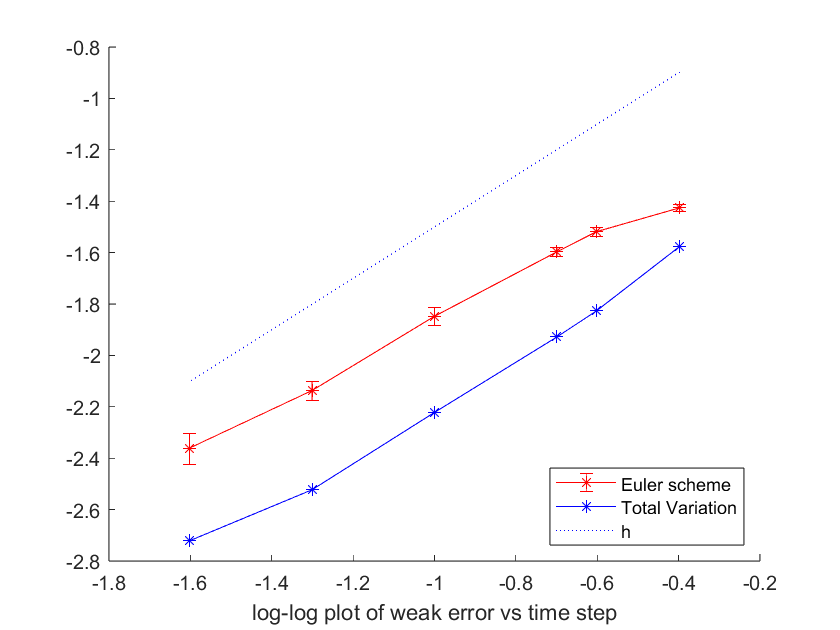}
\end{center}
\par
\captionof{figure}{Example 1: mixture of two univariate Gaussian distributions
(\ref{eq:uGm2}). The global error and the total variation distance for the
Euler method. Error bars indicate the Monte Carlo error. \label{fig:2Gauerr}}
\end{figure}

\begin{figure}[tbp]
\begin{center}
\par
\includegraphics[scale=0.5]{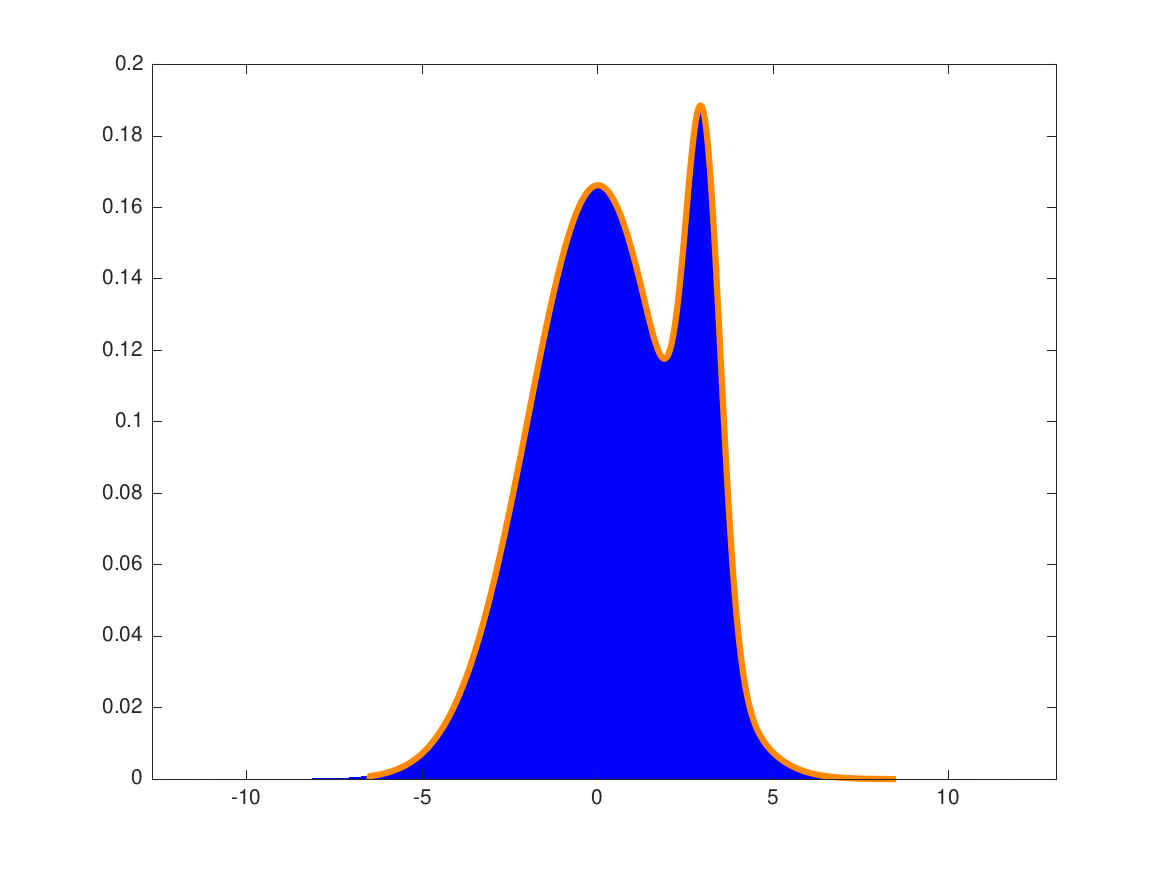}
\captionof{figure}{Example 1. Plot of the exact density $\rho (x)$ (\ref{eq:uGm2}) (orange) and the normalized histogram corresponding to the Euler
approximation (blue) with $h=0.025$, $T=100$ and $M=4 \cdot10^8$. \label{fig:2Gaudist}}
\end{center}
\end{figure}

The results are presented in Fig.~\ref{fig:2Gauerr} and Table~\ref{Texp1}.
The error shown is computed as $e:=\hat{\varphi}^{erg}-\varphi ^{erg}.$ We
also evaluate the total variation distance (TVD) between the exact density $%
\rho (x)$ and the density corresponding to the Euler approximation. We see
that both the weak-sense error for $\varphi (x)=x^{2}$ and the total
variation converge with order one in $h.$ We have not obtained a theoretical
result for convergence in total variation, which is an interesting topic for
future research. Figure~\ref{fig:2Gaudist} compares the exact density $\rho
(x)$ (\ref{eq:uGm2}) and the normalized histogram corresponding to the Euler
approximation.\medskip

\begin{table}[h] \centering%
\caption{Example 1: mixture of two univariate Gaussian distributions (\ref{eq:uGm2}).
The global error for the estimator $\hat \varphi$ of $\varphi^{erg}$ and the total variation.
\label{Texp1}}\smallskip $%
\begin{tabular}{lllllll}
\hline
$h$ & $M$ & $T$ & $\hat{\varphi}^{erg}$ & $\hat{\varphi}^{erg}-\varphi
^{erg} $ & $2\sqrt{\bar{D}_{M}/M}$ & TVD \\ \hline
$0.4$ & $10^{8}$ & $100$ & $4.9125$ & $0.0375$ & $0.0012$ & $0.0265$ \\
\hline
$0.25$ & $10^{8}$ & $100$ & $4.9053$ & $0.0303$ & $0.0012$ & $0.0149$ \\
\hline
$0.2$ & $10^{8}$ & $100$ & $4.9003$ & $0.0253$ & $0.0011$ & $0.0118$ \\
\hline
$0.1$ & $10^{8}$ & $100$ & $4.8892$ & $0.0142$ & $0.0011$ & $0.0060$ \\
\hline
$0.05$ & $4\cdot 10^{8}$ & $100$ & $4.8823$ & $0.0073$ & $0.0006$ & $0.0030$
\\ \hline
$0.025$ & $4\cdot 10^{8}$ & $100$ & $4.8794$ & $0.0044$ & $0.0006$ & $0.0019$
\\ \hline
\end{tabular}%
$%
\end{table}%

\begin{figure}[tbp]
\begin{center}
\par
\includegraphics[width=8cm]{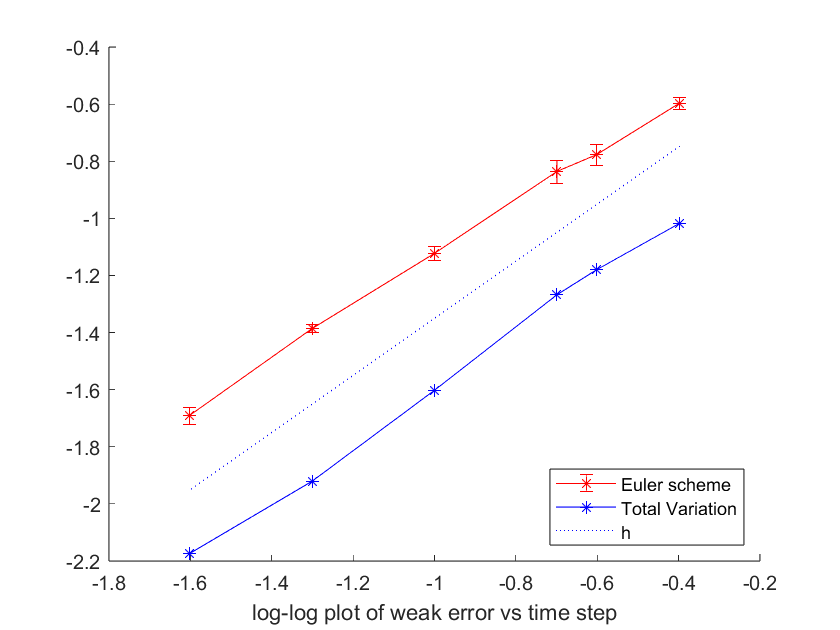}
\end{center}
\par
\captionof{figure}{Example 2: mixture distribution (\ref{eq:exa2}).
 Global error and the total variation distance for the
Euler method. Error bars indicate the Monte Carlo error. \label{fig:Example2Conv}}
\end{figure}

\noindent \textbf{Example 2.} Consider a mixture of two univariate Gaussian
distributions and a distribution corresponding to a nonconvex potential $%
U(x;3)=\beta (x^{4}-4x^{2})$:
\begin{eqnarray}
\rho (x)=&&\frac{1}{\mathcal{Z}}\left[ \alpha _{1}\exp \left( -\frac{(x-%
\mathrm{m}_{1})^{2}}{2\sigma _{1}^{2}}\right) \right .   \label{eq:exa2} \\
&&\left . +\alpha _{2}\exp \left( -\frac{%
(x-\mathrm{m}_{2})^{2}}{2\sigma _{2}^{2}}\right) +\alpha _{3}\exp \left(
-\beta (x^{4}-4x^{2})\right) \right] ,\ \ x\in \mathbb{R}, \notag
\end{eqnarray}%
where $\mathcal{Z}$ is the normalisation constant. The parameters used are
given in Table~\ref{Tparam2} and for test purposes we take $\varphi
(x;m)=x^{2}.$ The exact value $\varphi ^{erg}=6.98355$ $(5$ d.p.$).$ The
results are given in Table~\ref{Texp2} and Fig.~\ref{fig:Example2Conv} which
show first order convergence in $h$ of the Euler method (\ref{eq:hmeth}).

\begin{table}[h] \centering%
\caption{Example 2. The parameters of the mixture (\ref{eq:exa2}).
\label{Tparam2}}\smallskip
\begin{tabular}{l|lll}
\hline
$m$ & $1$ & $2$ & $3$ \\ \hline
$\alpha _{m}$ & $0.8$ & $1$ & $0.4$ \\ \hline
$c_{m}$ & $3.5$ & $-3$ &  \\ \hline
$\sigma _{m}$ & $1$ & $0.6$ &  \\ \hline
$\beta $ &  &  & $0.25$ \\ \hline
\end{tabular}%
\end{table}%

\begin{table}[h] \centering%
\caption{Example 2: mixture distribution (\ref{eq:exa2}).
The global error for the estimator $\hat \varphi$ of $\varphi^{erg}$ and the total variation.
\label{Texp2}}{}$%
\begin{tabular}{lllllll}
\hline
$h$ & $M$ & $T$ & $\hat{\varphi}^{erg}$ & $\hat{\varphi}^{erg}-\varphi ^{erg}
$ & $2\sqrt{\bar{D}_{M}/M}$ & TVD \\ \hline
$0.4$ & $10^{7}$ & $200$ & $6.731$ & $-0.253$ & $0.013$ & $0.096$ \\ \hline
$0.25$ & $10^{6}$ & $200$ & $6.816$ & $-0.167$ & $0.014$ & $0.066$ \\ \hline
$0.2$ & $10^{6}$ & $200$ & $6.837$ & $-0.146$ & $0.014$ & $0.054$ \\ \hline
$0.1$ & $10^{7}$ & $200$ & $6.9082$ & $-0.0754$ & $0.0043$ & $0.025$ \\
\hline
$0.05$ & $10^{8}$ & $200$ & $6.9424$ & $-0.0411$ & $0.0014$ & $0.012$ \\
\hline
$0.025$ & $10^{8}$ & $200$ & $6.9631$ & $-0.0204$ & $0.0014$ & $0.0067$ \\
\hline
\end{tabular}%
$%
\end{table}%

\begin{figure}[tbp]
\begin{center}
\par
\includegraphics[scale=0.5]{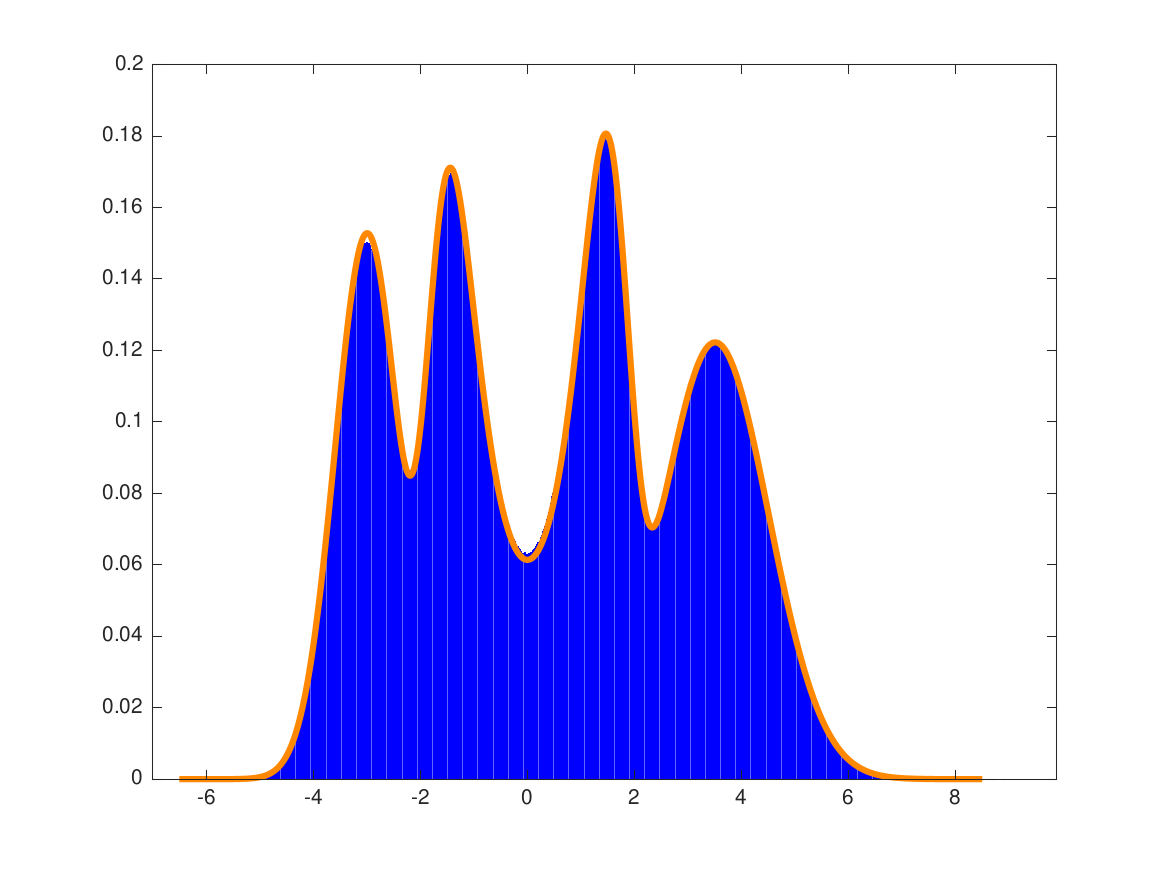}
\captionof{figure}{Example 2. Plot of the exact density $\rho (x)$ (\ref{eq:exa2}) (orange) and the normalized histogram corresponding to the Euler
approximation (blue) with $h=0.025$, $T=200$ and $M=10^8$. \label{fig:Example2dist}}
\end{center}
\end{figure}

The gradient of the potential $U(x;3)$ is growing faster than linearly at
infinity and hence Assumption~\ref{ass:erg2} is not satisfied. As it is well
known (see \cite{Stua,Tal99,MT} and references therein), explicit Euler-type
schemes are divergent when SDEs' coefficients do not satisfy the global
Lipschitz condition. The concept of rejecting exploding trajectories from
\cite{GNT04,MilTre07} allows us to apply any numerical method to SDEs with
nonglobally Lipschitz coefficients for estimating averages $E\varphi (X(T))$%
. The concept is based on rejecting trajectories which leave a sufficiently
large ball $S_{R}:=\{x:|x|<R\}$ during the time $T$ that contributes a very
small additional error to the overall simulation error thanks to probability
of $X(t),$ $t\in \lbrack 0,T],$ being exponentially small \cite%
{HAS,GNT04,MilTre07,MT}. Here we implement this concept as follows: if $%
|X_{k}|\geq 100$ for any $k=1,\ldots ,N,$ then we set $X_{n}=0$ for $n\geq k.
$ For $h=0.4$ there were $3.5\%$ rejected trajectories, for $h=0.25$ there
were just $3$ rejected trajectories out of $10^{6}$ and none for smaller
time steps. The experimental results confirm that the concept of rejecting
exploding trajectories can be successfully applied in the considered SDEwS
setting.  \medskip

\begin{figure}[tbp]
\begin{center}
\par
\includegraphics[width=11cm]{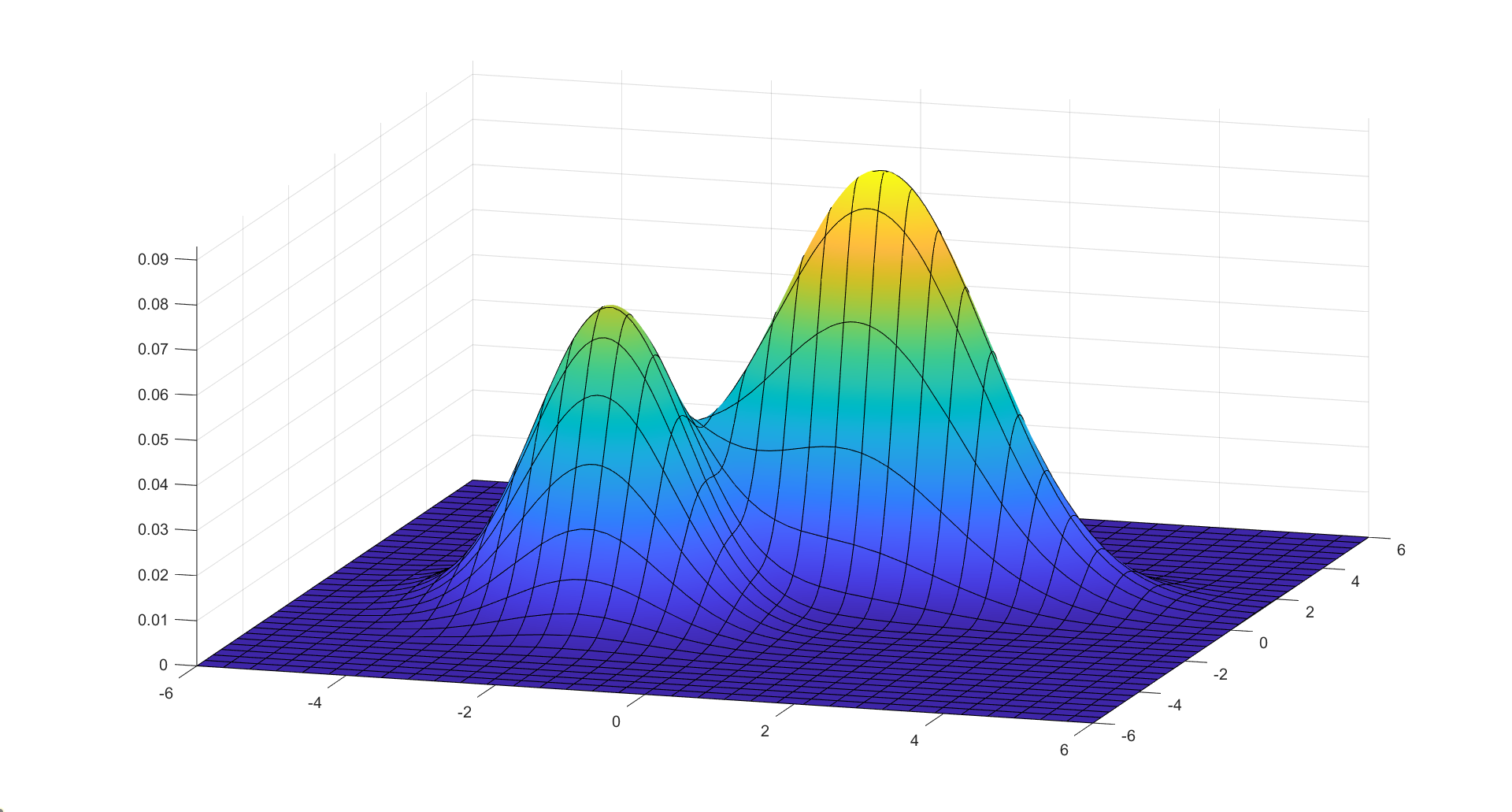}
\end{center}
\par
\captionof{figure}{Example 3: mixture distribution $\rho$  from (\ref{eq:exa3}) with parameters as in Table~\ref{Tparam3}.
  \label{fig:Example3rho}}
\end{figure}

\noindent \textbf{Example 3. }Consider a mixture of two multivariate
Gaussian distributions:
\begin{eqnarray}
\rho (x)=&&\frac{1}{\mathcal{Z}}\left[ \alpha _{1}\exp \left( -\frac{1}{2}(x-%
\mathrm{m}_{1})^{\top }\Sigma _{1}^{-1}(x-\mathrm{m}_{1})\right) \right . \label{eq:exa3} \\
&&\left . +\alpha
_{2}\exp \left( -\frac{1}{2}(x-\mathrm{m}_{2})^{\top }\Sigma _{2}^{-1}(x-%
\mathrm{m}_{2})\right) \right] ,\ \ x\in \mathbb{R}^{2},  \notag
\end{eqnarray}%
where $\mathcal{Z}$ is the normalisation constant, $\mathrm{m}_{i}$ are mean
vectors, and $\Sigma _{i}$ are covariance matrices:
\begin{equation*}
\Sigma _{i}=\left[
\begin{array}{cc}
\sigma _{1i} & r_{i} \\
r_{i} & \sigma _{2i}%
\end{array}%
\right] .
\end{equation*}%
The parameters used are given in Table~\ref{Tparam3} and for test purposes
we take $\varphi (x;m)=|x|^{2}.$ The density $\rho (x)$ is plotted in Fig.~%
\ref{fig:Example3rho}.

\begin{table}[h] \centering%
\caption{Example 3. The parameters of the mixture (\ref{eq:exa3}).
\label{Tparam3}}\smallskip
\begin{tabular}{l|ll}
\hline
$i$ & $1$ & $2$ \\ \hline
$\alpha _{i}$ & $0.7$ & $0.5$ \\ \hline
$\mathrm{m}_{i}$ & $(1,1)^{\top }$ & $(-2,-1)^{\top }$ \\ \hline
$\sigma _{1i}$ & $2$ & $1$ \\ \hline
$\sigma _{2i}$ & $0.5$ & $1$ \\ \hline
$r_{i}$ & $0.1$ & $-0.1$ \\ \hline
\end{tabular}%
\end{table}%

\begin{figure}[tbp]
\begin{center}
\par
\includegraphics[width=12cm]{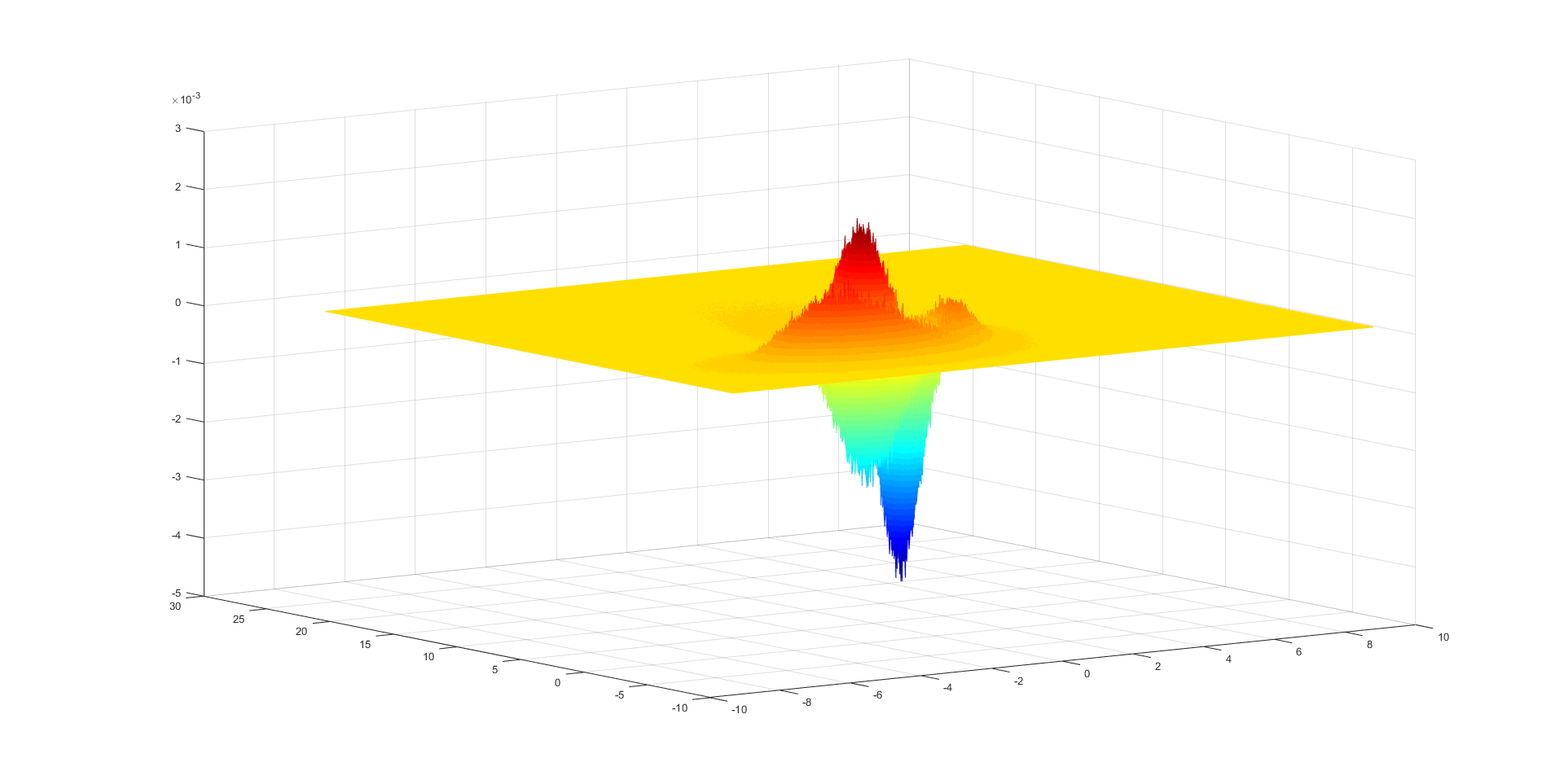}
\end{center}
\par
\captionof{figure}{Example 3: The difference between the normalized histogram corresponding
to the Euler approximation with $h=0.1$, $M=10^9$, $T=200$ and the exact density $\rho$  from (\ref{eq:exa3}) with parameters as in Table~\ref{Tparam3}.
  \label{fig:3Error3D}}
\end{figure}

The exact value $\varphi ^{erg}=5.541667$ $(6$ d.p.$).$ The results of the
experiments are presented in Table~\ref{Texp3} and Fig.~\ref%
{fig:Example3Conv}, which again show first order convergence in $h$ of the
Euler method (\ref{eq:hmeth}) as expected. Figure~\ref{fig:3Error3D} gives
the difference between the normalized histogram corresponding to the Euler
approximation and the exact density $\rho (x)$ (\ref{eq:exa3}).

\begin{table}[h] \centering%
\caption{Example 3: mixture of two univariate Gaussian distributions (\ref{eq:uGm2}).
The global error for the estimator $\hat \varphi$ of $\varphi^{erg}$ and the total variation.
\label{Texp3}}\smallskip $%
\begin{tabular}{lllllll}
\hline
$h$ & $M$ & $T$ & $\hat{\varphi}^{erg}$ & $\hat{\varphi}^{erg}-\varphi
^{erg} $ & $2\sqrt{\bar{D}_{M}/M}$ & TVD \\ \hline
$0.5$ & $10^{6}$ & $200$ & $5.8559$ & $0.3142$ & $0.0101$ & $0.0744$ \\
\hline
$0.4$ & $10^{6}$ & $200$ & $5.7798$ & $0.2381$ & $0.0099$ & $0.0581$ \\
\hline
$0.32$ & $10^{7}$ & $200$ & $5.7279$ & $0.1862$ & $0.0031$ & $0.0430$ \\
\hline
$0.25$ & $10^{8}$ & $200$ & $5.6808$ & $0.1392$ & $0.0010$ & $0.0302$ \\
\hline
$0.2$ & $4\cdot 10^{8}$ & $200$ & $5.6528$ & $0.1112$ & $0.0005$ & $0.0238$
\\ \hline
$0.16$ & $4\cdot 10^{8}$ & $200$ & $5.6292$ & $0.0876$ & $0.0005$ & $0.0197$
\\ \hline
$0.1$ & $10^{9}$ & $200$ & $5.5961$ & $0.0544$ & $0.0003$ & 0.0131 \\ \hline
\end{tabular}%
$%
\end{table}%

\begin{figure}[tbp]
\begin{center}
\par
\includegraphics[width=8cm]{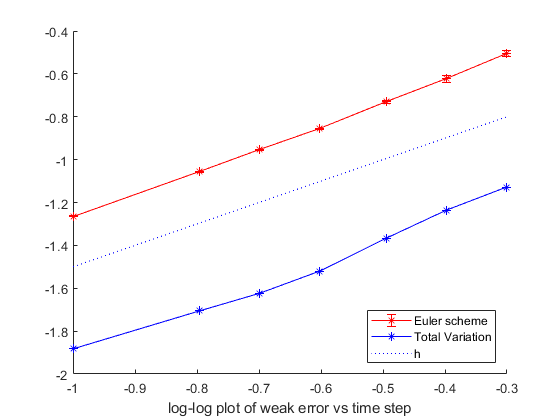}
\end{center}
\par
\captionof{figure}{Example 3: mixture distribution (\ref{eq:exa3}).
 Global error and the total variation distance for the
Euler method. Error bars indicate the Monte Carlo error. \label{fig:Example3Conv}}
\end{figure}

\section{Proofs\label{sec:proof}}

The proofs of both Theorems~\ref{thm:ft} and \ref{thm:erg1} follow the same
path: we first prove a second-order weak error bound for the corresponding
one-step approximation (Lemmas~\ref{lem:onestft} and~\ref{lem:onesterg}) and
then appropriately re-write the global error to sum up the local errors.

\subsection{Finite-time convergence\label{sec:ftconv}}

The proof of the next lemma on boundedness of moments of the Euler scheme
under weaker conditions than Assumption~\ref{ass:ft} can be found e.g. in
\cite[Lemma 3.1]{SDEwMS7}.

\begin{lemma}
\label{lem:momft}Let Assumption~\ref{ass:ft} hold. Then the method (\ref%
{eq:Eul}) has finite moments, i.e., for all $p\geq 1:$
\begin{equation*}
E\left( |X_{k}|^{2p}\right) \leq K(1+|x|^{2p}),
\end{equation*}%
where $K>0$ does not depend on $h$ and $x.$
\end{lemma}

The following error gives the one-step error estimate.

\begin{lemma}
\label{lem:onestft}Let Assumption~\ref{ass:ft} hold. Then the one-step error
of the method (\ref{eq:Eul}) is estimated as
\begin{equation*}
\left\vert Eu(t+h,X_{1};\mu _{1})-u(t,x;m)\right\vert \leq
C(1+|x|^{\varkappa })h^{2},
\end{equation*}%
where $u(t+h,x,m),$ $m\in \mathcal{M},$ is a solution of (\ref{eq:BKE}), $%
C>0 $ and $\varkappa \geq 1$ do not depend on $h,$ $x$ and $m.$
\end{lemma}

\noindent \textbf{Proof.} Using Taylor's expansion, the conditions (\ref{eq:rn}), and the backward Kolmogorov
equation (\ref{eq:BKE}), we obtain
\begin{eqnarray}
&&Eu(t+h,X_{1};\mu _{1})-u(t,x;m)  \label{eq:lem11} \\
&=&h\sum_{j\neq m}q_{mj}(x)Eu(t+h,X_{1};j)+\left( 1-hq_{m}(x)\right)
Eu(t+h,X_{1};m)-u(t,x;m)  \notag \\
&=&h\sum_{j\neq
m}q_{mj}(x)u(t,x;j)-hq_{m}(x)u(t,x;m)+Eu(t+h,X_{1};m) \notag \\
&&-u(t,x;m)+r_{1}(t,x)  \notag \\
&=&h\rf{{\cal L}u(t,x;m)}  +r_{2}(t,x) = r_{2}(t,x),  \notag
\end{eqnarray}%
where \rf{$r_1$ is a combination of first order derivative in time and first and second order derivatives in space of $u$ multiplied by corresponding $q_{ij}$ and the coefficients of (\ref{eq:sde}); $r_2$ in addition contains a combination of
$\frac{\partial^{l+|i|} u}{\partial t^{l} \partial x^{i_1} \cdots \partial x^{i_d}}$ with $l=0$, $|i|=2,3,4$; $l=1$, $|i|=1,2$; $l=2$, $|i|=0$ appropriately multiplied by the coefficients of (\ref{eq:sde}); and, thanks to Assumption~\ref{ass:ft} and the estimate (\ref{eq:FKest}), the reminders are estimated as}
\begin{equation*}
|r_{i}(t,x)|\leq C(1+|x|^{\varkappa })h^{2}.
\end{equation*}%
$\ \square \medskip $

\textbf{Proof of Theorem~\ref{thm:ft}}. We re-write the global error using
the Feynman-Kac formula (\ref{eq:Kac}):%
\begin{gather*}
\left\vert E\varphi (X_{N};\mu _{N})-E\varphi (X_{t,x,m}(T);\mu
_{t,x,m}(T))\right\vert =\left\vert Eu(t_{N},X_{N};\mu
_{N})-u(t,x;m)\right\vert \\
=\left\vert E\sum_{k=0}^{N-1}(u(t_{k+1},X_{k+1};\mu
_{k+1})-u(t_{k},X_{k};\mu _{k}))\right\vert \\
=\left\vert \sum_{k=0}^{N-1}E\left( E\left[ u(t_{k+1},X_{k+1};\mu
_{k+1})-u(t_{k},X_{k};\mu _{k})|X_{k};\mu _{k}\right] \right) \right\vert \\
\leq \sum_{k=0}^{N-1}E\left\vert E\left[ u(t_{k+1},X_{k+1};\mu
_{k+1})-u(t_{k},X_{k};\mu _{k})|X_{k};\mu _{k}\right] \right\vert \leq
Ch^{2}\sum_{k=0}^{N-1}E(1+|X_{k+1}|^{\varkappa }) \\
\leq Ch(1+|x|^{\varkappa }),
\end{gather*}%
where we applied Lemmas~\ref{lem:momft}\ and~\ref{lem:onestft}. $\ \square $

\subsection{Ergodic limits convergence\label{sec:ergproof}}

We start with proving an estimate for the moments of the method (\ref%
{eq:hmeth}).

\begin{lemma}
\label{lem:momerg}Let \rf{Assumption~\ref{ass:erg2}} hold. Then the method (\ref%
{eq:hmeth}) has finite moments, i.e., for sufficiently small $h>0$ and all $%
p\geq 1:$
\begin{equation*}
E\left( |X_{k}|^{2p}\right) \leq K(1+|x|^{2p}e^{-\lambda t_{k}}),
\end{equation*}%
where $K>0$ does not depend on $h,$ $T$ and $x.$
\end{lemma}

\noindent \textbf{Proof. }We note that the constant $K>0$ is changing from
line to line. We have%
\begin{eqnarray*}
E|X_{k+1}|^{2p} & \rf{\leq} &E|X_{k}|^{2p}+pE|X_{k}|^{2p-2}[2(X_{k}\cdot
(X_{k+1}-X_{k}))+(2p-1)|X_{k+1}-X_{k}|^{2}] \\
&&+K\sum_{l=3}^{2p}E|X_{k}|^{2p-l}|X_{k+1}-X_{k}|^{l},
\end{eqnarray*}%
where $K>0$ depends on $p$ only. Using that $|\nabla U(x;m)|\leq K(1+|x|)$, 
(\ref{eq:rn}) and (\ref{eq:condU}), we obtain
\begin{eqnarray*}
E[X_{k}\cdot (X_{k+1}-X_{k})|X_{k},\mu _{k}] &\leq &-\frac{h}{2}X_{k}\cdot
\nabla U(X_{k},\mu _{k})\leq -\frac{h}{2}c_{1}|X_{k}|^{2}+\frac{h}{2}c_{2},\,
\\
E[|X_{k+1}-X_{k}|^{2}|X_{k},\mu _{k}] &\leq &Kh^{2}|X_{k}|^{2}+h, \\
E[|X_{k+1}-X_{k}|^{l}|X_{k},\mu _{k}] &\leq
&Kh^{l}|X_{k}|^{l}+Kh^{l/2},\;l\geq 3,
\end{eqnarray*}%
with $K>0$ independent of $h$.

Then, using Young's inequality (i.e., for $y,z\geq 0:$ $yz\leq \epsilon ^{i}%
\frac{y^{i}}{i}+\epsilon ^{-j}\frac{z^{j}}{j}$ with $\epsilon >0$, $i,j>1$
and $\frac{1}{i}+\frac{1}{j}=1)$, we get (where needed with a new $K>0$
independent of $h):$%
\begin{gather*}
2pE|X_{k}|^{2p-2}(X_{k}\cdot (X_{k+1}-X_{k}))\leq
-pc_{1}hE|X_{k}|^{2p}+phKE|X_{k}|^{2p-2} \\
\leq -pc_{1}hE|X_{k}|^{2p}+(p-1)c_{1}h\mathbb{E}|X_{k}|^{2p}+hK^{p}\left(
\frac{1}{c_{1}}\right) ^{p-1} \\
=-c_{1}hE|X_{k}|^{2p}+hK^{p}\left( \frac{1}{c_{1}}\right) ^{p-1},
\end{gather*}%
\begin{gather*}
p(2p-1)E|X_{k}|^{2p-2}|X_{k+1}-X_{k}|^{2}\leq
p(2p-1)Kh^{2}|X_{k}|^{2p}+hp(2p-1)E|X_{k}|^{2p-2} \\
\leq Kh^{2}|X_{k}|^{2p}+\frac{c_{1}}{2}hE|X_{k}|^{2p}+\frac{[(2p-1)(p-1)]^{p}%
}{p-1}\left( \frac{2}{c_{1}}\right) ^{p-1}h,
\end{gather*}%
and for $l\geq 3$%
\begin{gather*}
E|X_{k}|^{2p-l}|X_{k+1}-X_{k}|^{l}\leq
Kh^{l}|X_{k}|^{2p}+Kh^{l/2}E|X_{k}|^{2p-l} \\
\leq Kh^{l}|X_{k}|^{2p}+\frac{c_{1}h^{l/2}}{2p}\mathbb{E}%
|X_{k}|^{2p}+Kh^{l/2}
\end{gather*}%
Thus
\begin{equation*}
E|X_{k+1}|^{2p}\leq E|X_{k}|^{2p}-\frac{c_{1}}{2}h(1-\sqrt{h}%
-Kh)E|X_{k}|^{2p}+Kh,
\end{equation*}%
whence by Gronwald's lemma%
\begin{gather*}
E|X_{k+1}|^{2p}\leq \left( 1-\frac{c_{1}h}{2}(1-\sqrt{h}-Kh)\right)
^{k}|x|^{2p}\\
+\frac{2K}{c_{1}(1-\sqrt{h}-Kh)}\left[ 1-\left( 1-\frac{c_{1}h}{2%
}(1-\sqrt{h}-Kh)\right) ^{k}\right] .
\end{gather*}

$\ \square \medskip $

The local error of the method (\ref{eq:hmeth}) is given in the next lemma.

\begin{lemma}
\label{lem:onesterg}Let Assumption~\ref{ass:erg2} hold. Then the one-step
error of the method (\ref{eq:hmeth}) is estimated as
\begin{equation*}
\left\vert Eu(t+h,X_{1};\mu _{1})-u(t,x;m)\right\vert \leq
C(1+|x|^{\varkappa })e^{-\lambda (T-t)}h^{2},
\end{equation*}%
where $u(t+h,x,m),$ $m\in \mathcal{M},$ is a solution of (\ref{eq:BKE2}), $%
C>0$ and $\varkappa \geq 1$ do not depend on $T$ and $h.$
\end{lemma}

\noindent \textbf{Proof. }Analogously to (\ref{eq:lem11}) but using (\ref%
{eq:BKE2}), we have
\begin{equation*}
Eu(t+h,X_{1};\mu _{1})-u(t,x;m)=r_{3}(t,x),
\end{equation*}%
where the reminder $r_{3}(t,x)$ is, \rf{analogously to $r_2$ in (\ref{eq:lem11}),} a sum of terms containing derivatives of $%
u$ and the coefficients of the SDEwS (\ref{eq:erg}), (\ref{eq:erg3}) and
hence by (\ref{eq:FKerg})%
\begin{equation*}
\left\vert r_{3}(t,x)\right\vert \leq C(1+|x|^{\varkappa })e^{-\lambda
(T-t)}h^{2},
\end{equation*}%
where $C>0\ $is independent of $h,$ $t,$ $T$ and $x$. \ $\square \medskip $

\textbf{Proof of Theorem~\ref{thm:erg1}}. Using (\ref{PA34}) and Lemmas~\ref%
{lem:onesterg} and \ref{lem:momerg}, we obtain
\begin{eqnarray*}
&&\left\vert E\varphi (X_{N};\mu _{N})-\varphi ^{erg}\right\vert \\
&\leq &\left\vert E\varphi (X_{N};\mu _{N})-E\varphi (X_{0,x,m}(T);\mu
_{0,x,m}(T))\right\vert +\left\vert E\varphi (X_{0,x,m}(T);\mu
_{0,x,m}(T))-\varphi ^{erg}\right\vert \\
&\leq &\sum_{k=0}^{N-1}E\left\vert E\left[ u(t_{k+1},X_{k+1};\mu
_{k+1})-u(t_{k},X_{k};\mu _{k})|X_{k};\mu _{k}\right] \right\vert +C\left(
1+|x|^{\kappa }\right) e^{-\lambda T} \\
&\leq &Ch^{2}\sum_{k=0}^{N-1}e^{-\lambda (T-t_{k})}(1+E|X_{k}|^{\varkappa
})+C\left( 1+|x|^{\kappa }\right) e^{-\lambda T}\leq C\left[
h+(1+|x|^{\kappa })\exp (-\lambda T)\right] .
\end{eqnarray*}%
$\ \square \medskip $

\textbf{Proof of Theorem~\ref{thm:erg2}}. Thanks to (\ref{eq:FKerg}), we get
\begin{eqnarray}
&&\left\vert h\sum_{l=1}^{L}E\varphi (X_{0,x,m}(t_{l});\mu
_{0,x,m}(t_{l}))-\int\limits_{0}^{\tilde{T}}E\varphi (X_{0,x,m}(s);\mu
_{0,x,m}(s))ds\right\vert  \label{eq:prThm3_1} \\
&=&\left\vert h\sum_{l=1}^{L}u(0,x;m;t_{l})-\int\limits_{0}^{\tilde{T}%
}u(0,x;m;s)ds\right\vert \notag \\
&&=\left\vert
h\sum_{l=1}^{L}u(-t_{l},x;m;0)-\int\limits_{0}^{\tilde{T}}u(-s,x;m;0)ds%
\right\vert  \notag \\
&\leq &Ch^{2}\sum_{l=1}^{L}\left( 1+|x|^{\kappa }\right) e^{-\lambda
(T-t_{l})}\leq Ch\left( 1+|x|^{\kappa }\right) .  \notag
\end{eqnarray}

Using Theorem~\ref{thm:erg1}, (\ref{eq:prThm3_1}) and\textbf{\ }(\ref%
{eq:biasTime}), we obtain%
\begin{eqnarray*}
&&\left\vert E\check{\varphi}^{erg}-\varphi ^{erg}\right\vert \leq
\left\vert \frac{1}{L}\sum_{l=1}^{L}\left[ E\varphi (X_{l};\mu
_{l})-E\varphi (X_{0,x,m}(t_{l});\mu _{0,x,m}(t_{l}))\right] \right\vert \\
&&+\left\vert \frac{1}{L}\sum_{l=1}^{L}E\varphi (X_{0,x,m}(t_{l});\mu
_{0,x,m}(t_{l}))-\frac{1}{\tilde{T}}E\int\limits_{0}^{\tilde{T}}\varphi
(X_{0,x,m}(s);\mu _{0,x,m}(s))ds\right\vert \\
&&+\left\vert \frac{1}{\tilde{T}}E\int\limits_{0}^{\tilde{T}}\varphi
(X_{0,x,m}(s);\mu _{0,x,m}(s))ds-\varphi ^{erg}\right\vert \\
&\leq &\frac{1}{L}\sum_{l=1}^{L}C\left[ h+(1+|x|^{\kappa })\exp (-\lambda
t_{l})\right] +C\frac{h}{\tilde{T}}\left( 1+|x|^{\kappa }\right) +C\left(
1+|x|^{\kappa }\right) \frac{1}{\tilde{T}} \\
&\leq &C\left( h+\frac{1+|x|^{\kappa }}{\tilde{T}}\right) .
\end{eqnarray*}%
\ $\square $

\begin{remark}\label{rem:Poisson}
The following system of Poisson equations is associated with SDEwS (\ref%
{eq:erg}), (\ref{eq:erg3}) with $Q(x)$ satisfying (\ref{eq:q_choice}):
\begin{equation}
\mathcal{L}u(x;m)=\varphi (x;m)-\varphi _{m}^{erg},\ x\in \mathbb{R}^{d},\
m\in \mathcal{M},  \label{eq:Poi}
\end{equation}%
where $\mathcal{L}$ is as in (\ref{eq:BKE2}) and
\begin{equation*}
\varphi _{m}^{erg}=\int_{\mathbb{R}^{d}}\varphi (x;m)\rho (x;m)dx
\end{equation*}%
with $\rho (x;m)$ from (\ref{eq:mix}). If $U(x,m)\in C^{3}(\mathbb{R}^{d}),$
$m\in \mathcal{M},$ and satisfy (\ref{eq:condU}), and the elements $q_{ij}(x)
$ of the matrix $Q(x)$ are positive functions bounded in $\mathbb{R}^{d}$
and belong to $C^{2}(\mathbb{R}^{d}),$ and the functions $\varphi (x;m),$ $%
m\in \mathcal{M},$ and their derivatives $\frac{\partial ^{|j|}\varphi }{%
\partial x^{j_{1}}\cdots \partial x^{j_{d}}},$ $|j|\leq 2,$ are continuous
and belong to the class $\mathbf{F,}$ then the solution to the system of
elliptic PDEs (\ref{eq:Poi}) satisfies the following estimate for some $%
\kappa \geq 1$ \rf{\cite[Chap. 9, Sec. 7]{Friedman}}:
\begin{equation*}
\sum\limits_{|j|=0}^{4}\left\vert \frac{\partial ^{|j|}}{\partial
x^{j_{1}},\dots ,\partial x^{j_{d}}}u(x;m)\right\vert \leq K(1+|x|^{\kappa
}),
\end{equation*}%
where $K>0$ does not depend on $x$. The Poisson equation can be
used in proving properties of the time-averaging estimator $\check{\varphi}%
^{erg}$ analogously to \cite{MST10} (see also \cite{LST23,manifold}).
\end{remark}

\section*{Acknowledgment}
For the purpose of open access, the author  applied a
Creative Commons Attribution (CC-BY) license to any Author Accepted Manuscript version arising.


\begin{thebibliography}{99}
\bibitem{SDEwMS3} J. Bao, J. Shao, C. Yuan. Approximation of invariant
measures for regime-switching diffusions. \textit{Potential Analysis}
\textbf{44} (2016), 707--727.

\bibitem{manifold} K. Bharath, A. Lewis, A. Sharma, M.V. Tretyakov. Sampling
and estimation on manifolds using the Langevin diffusion. arXiv:2312.14882,
2023.

\bibitem{mix2} N. Bouguila, W. Fan. \textit{Mixture models and applications}. Springer, 2020.

\rf{
\bibitem{CMH} B. Cloez, M. Hairer. Exponential ergodicity for Markov processes with random switching.
\textit{Bernoulli} \textbf{21} (2015), 505--536.
}

\bibitem{mix3} B.S. Everitt, D.J. Hand. \textit{Finite mixture distributions}%
. Springer, 1981.

\bibitem{Friedman} A. Fridman. \textit{Partial differential equations of
parabolic type}. Prentice-Hall, 1964.

\bibitem{mix1} S. Fruehwirth-Schnatter. \textit{Finite mixture and Markov
switching models}. Springer, 2007.

\bibitem{Lot24} Z. Gou, X. Tu, S.V. Lototsky, R. Ghanem. Switching
diffusions for multiscale uncertainty quantification. \textit{Inter. J.
Non-Linear Mechanics} (2024), 104793.

\bibitem{HJ15} M. Hutzenthaler and A. Jentzen. \textit{Numerical
approximations of stochastic differential equations with non-globally
Lipschitz continuous coefficients.} Mem. Amer. Math. Soc., vol. 236. AMS,
Providence, 2015.

\bibitem{HAS} R.Z. Khasminskii. Stochastic stability of differential
equations. Springer, 2012.

\bibitem{LSU88} O.A. Ladyzhenskaya, V.A. Solonnikov, N.N. Ural'ceva, \textit{%
Linear and quasilinear equations of parabolic type}. Amer. Math.
Soc., 1988.

\bibitem{LST23} B. Leimkuhler, A. Sharma, M. V. Tretyakov. Simplest random
walk for approximating Robin boundary value problems and ergodic limits of
reflected diffusions. \textit{Ann. Appl. Probab.} \textbf{33} (2023),
1904--1960.

\bibitem{SDEwMS4} X. Li, Q. Ma, H. Yang, C. Yuan. The numerical invariant
measure of stochastic differential equations with Markovian switching.
\textit{SIAM J. Numer. Anal.} \textbf{56} (2018), 1435--1455.

\bibitem{Maobook} X. Mao, C. Yuan. \textit{Stochastic differential equations
with Markovian switching}. Imperial College Press, 2006.

\bibitem{SDEwMS1} X. Mao, C. Yuan, G. Yin. Numerical method for stationary
distribution of stochastic differential equations with Markovian switching.
\textit{J. Comp. Appl. Math.} \textbf{174} (2005), 1--27.

\bibitem{Stua} J.C. Mattingly, A.M. Stuart, D.J. Higham. Ergodicity for SDEs
and approximations: Locally Lipschitz vector fields and degenerate noise.
\textit{Stoch. Proc. Appl.}, \textbf{101} (2002), 185--232.

\bibitem{MST10} J.C. Mattingly, A.M. Stuart, M.V. Tretyakov. Convergence of
numerical time-averaging and stationary measures via Poisson equations.
\textit{SIAM J. Numer. Anal.} \textbf{48} (2010), 552--577.

\bibitem{SDEwMS5} H. Mei, G. Yin. Convergence and convergence rates for
approximating ergodic means of functions of solutions to stochastic
differential equations with Markov switching. \textit{Stoch. Proc. Appl.}
\textbf{125} (2015), 3104--3125.

\bibitem{GN72} G.N. Milstein. Interaction of Markov processes. \textit{%
Theory Prob. Appl.} \textbf{17} (1972), 36--45.

\bibitem{GN78a} G.N. Milstein. A method with second order accuracy for the
integration of stochastic differential equations. \textit{Theor. Prob. Appl.}
\textbf{23} (1978), 414--419.

\bibitem{GN78b} G.N. Milstein. Probabilistic solution of linear systems of
elliptic and parabolic equations. \textit{Theor. Prob. Appl.} \textbf{23}
(1978), 851--855.

\bibitem{GN85} G.N. Milstein. Weak approximation of solutions of systems of
stochastic differential equations. \textit{Theor. Prob. Appl.} \textbf{30}
(1985), 706--721.

\bibitem{MT} G.N. Milstein, M.V. Tretyakov. \textit{Stochastic numerics for
mathematical physics. }2nd edition\textit{.} Springer, 2021.

\bibitem{GNT04} G.N. Milstein, M.V. Tretyakov. Numerical integration of
stochastic differential equations with nonglobally Lipschitz coefficients.
\textit{SIAM J. Numer. Anal.} \textbf{43\ (}2005), 1139--1154.

\bibitem{MilTre07} G.N. Milstein, M.V. Tretyakov. Computing ergodic limits
for Langevin equations. \textit{Phys. D} \textbf{229} (2007), 81--95.

\bibitem{RT96} G.O. Roberts, R.L. Tweedie. Exponential convergence of
Langevin distributions and their discrete approximations. \textit{Bernoulli}
\textbf{2} (1996), 341--363.

\bibitem{Skoroh} A.V. Skorohod. Asymptotic methods in the theory of
stochastic differential equations. AMS, 1989.

\bibitem{SDEwMS8} J. Shao. Invariant measures and Euler--Maruyama's
approximations of state-dependent regime-switching diffusions. \textit{SIAM
J. Control Optim.} \textbf{56} (2018), 3215--3238.

\bibitem{Tal86} D. Talay. Discr\'{e}tisation d'une \'{e}quation diff\'{e}%
rentielle stochastique et calcul approch\'{e} d'esp\'{e}rances de
fonctionnelles de la solution. \textit{Math. Model. Numer. Anal. (ESAIM)}
\textbf{20 }(1986), 141--179.

\bibitem{TT90} D. Talay, L. Tubaro. Expansion of the global error for
numerical schemes solving stochastic differential equations. \textit{Stoch.
Anal. Appl.} \textbf{8 }(1990), 483--509.

\bibitem{Tal90} D. Talay. Second-order discretization schemes for stochastic
differential systems for the computation of the invariant law. \textit{%
Stoch. Stoch. Reports}, \textbf{29 }(1990), 13--36.

\bibitem{Tal99} D. Talay. Stochastic Hamiltonian systems: exponential
convergence to the invariant measure, and discretization by the implicit
Euler scheme. \textit{Markov Proc. Relat. Fields}, \textbf{8} (2002),
163--198.

\rf{
\bibitem{Majda} X.T. Tong, A.J. Majda.
Moment bounds and geometric ergodicity of diffusions with random switching and unbounded transition rates.
\textit{Res. Math. Sci.} \textbf{3} (2016), Article 41.
}

\bibitem{Handy13} M.V. Tretyakov, Z. Zhang. A fundamental mean-square
convergence theorem for SDEs with locally Lipschitz coefficients and its
applications. \textit{SIAM J. Numer. Anal.} \textbf{51} (2013), 3135--3162.

\rf{
\bibitem{XuZhu}F. Xi, C. Zhu. On Feller and strong Feller properties and exponential ergodicity
of regime-switching jump diffusion processes with countable regimes.
\textit{SIAM J. Control Optimiz.} \textbf{55} (2017), 1789--1818.
}

\bibitem{SDEwMS7} G. Yin, X. Mao, C. Yuan, D. Cao. Approximation methods for
hybrid diffusion systems with state-dependent switching processes: numerical
algorithms and existence and uniqueness of solutions. \textit{SIAM J. Math.
Anal.} \textbf{41} (2010), 2335--2352.

\bibitem{Yinbook} G. Yin, C. Zhu. \textit{Hybrid switching diffusions:
properties and applications}. Springer, 2009.

\bibitem{SDEwMS2} G. Yin, C. Zhu. Properties of solutions of stochastic
differential equations with continuous-state-dependent switching. \textit{J.
Diff. Equations} \textbf{249} (2010), 2409--2439.
\end{thebibliography}
\end{document}